\documentclass[12pt]{article}
\usepackage{mathptmx}
\usepackage{amsfonts}
\newtheorem{Theorem}{Theorem}[section]

\newtheorem{Lemma}{Lemma}[section]

\newtheorem{Remark}{Remark}[section]

\title{Approximation by generalized Kantorovich sampling type series}
\author{A. Sathish Kumar\thanks{Department of Mathematics, Visvesvaraya National Institute of Technology, Nagpur, Nagpur-440010, India. \newline E-mail: mathsatish9@gmail.com}
 \and
 P. Devaraj\thanks{Department of Mathematics, Indian Institute of Science Education and Research, Thiruvananthapuram \newline E-mail: devarajp@iisertvm.ac.in}
  }
\date{}

\begin{document}
\maketitle
\bibliographystyle{plain}
\abstract{In the present article, we analyse the behaviour of a new family of Kantorovich type sampling operators $(K_w^{\varphi}f)_{w>0}.$ First, we give a Voronovskaya type theorem for these Kantorovich generalized sampling series and a corresponding quantitative version in terms of the first order of modulus of continuity.  Further, we study the order of approximation in $C({\mathbb{R}})$ (the set of all uniformly continuous and bounded functions on ${\mathbb{R}}$) for the family $(K_w^{\varphi}f)_{w>0}.$ Finally, we give some examples of kernels such as B-spline kernels and Blackman-Harris kernel to which the theory can be applied.

\endabstract



\section{Introduction}

The theory of generalized sampling series was first initiated by P. L. Butzer and his school \cite{Stens1} and \cite{Stens2}. In recent years, it is an attractive topic in approximation theory due to its wide range of applications, especially in signal and image processing. For $w>0,$ a generalized sampling series of a function $f:{\mathbb{R}}\rightarrow{\mathbb{R}}$ is defined by
\begin{eqnarray*}
(T_w^{\varphi}f)(x)&=&\sum_{k=-\infty}^{\infty}\varphi(wx-k)f\bigg(\frac{k}{w}\bigg), \,\ x\in{\mathbb{R}},
\end{eqnarray*}
where $\varphi$ is a kernel function on ${\mathbb{R}}.$ These type of operators have been studied by many authors ( e.g. \cite{a1}, \cite{a3}, \cite{a2},  \cite{a4}, \cite{a6}, \cite{a5},  \cite{a7}, \cite{a8}).

The Kantorovich type generalizations of approximation operators is an important subject in
approximation theory and they are the method to approximate Lebesgue integrable functions.
In the last few decades, the Kantorovich modifications of several operators were constructed and their approximation behavior studied, we mention some of the work in this direction e.g., \cite{gupta1,maria,PNA, Tuncer} etc.

In \cite{BVBS}, the authors have introduced the sampling Kantorovich operators and studied their rate of convergence in the general settings of Orlicz spaces. After that, Costarelli and Vinti \cite{COS} extended their study in the multivariate setting and obtained the rate of convergence for functions in Orlicz spaces. Danilo and Vinti \cite{Dani}, obtained the rate of approximation for the family of sampling Kantorovich operators in the uniform norm, for bounded uniformly continuous  functions belonging to Lipschitz classes and for functions in Orlicz spaces. Also, the nonlinear version of sampling Kantorovich operators has been studied in \cite{Dan} and
\cite{Zam}.

Altomare and Leonessa \cite{Alt} considered a new sequence of positive linear operators acting on the space of Lebesgue-integrable functions on the unit interval. Such operators include the Kantorovich operators as a particular case. Later, in order to obtain an approximation process for spaces of locally integrable functions on unbounded intervals, Altomare et. al. introduced and studied the generalized  Sz\'{a}sz-Mirakjan-Kantorovich operators in \cite{MMC}. Also in \cite{Cap}, the authors obtained some qualitative properties and an asymptotic formula for such a sequence of operators.

We consider the generalized Kantorovich type sampling series.
Let $\{a_k\}_{k\in {\mathbb{Z}}}$ and $\{b_k\}_{k\in {\mathbb{Z}}}$ be two sequences of real numbers such that $b_k-a_k=\Delta_{k}>0$ for $k\in {\mathbb{Z}}$. In this paper, we analyse the approximation properties of the following type of generalized Kantorovich sampling series
\begin{eqnarray}\label{e1}
(K_w^{\varphi}f)(x)&=&\sum_{k=-\infty}^{\infty}\frac{w}{\Delta_{k}}\varphi(wx-k)\int_{\frac{k+a_k}{w}}^{\frac{k+b_k}{w}}f(u)du,
\end{eqnarray}
$f\in C({\mathbb{R}})$ (The class of all uniformly continuous and bounded functions on ${\mathbb{R}}$).

In the present paper, first we obtain asymptotic formula and their quantitative estimate for the operators  $(K_w^{\varphi}f)_{w>0}.$ Further, the order of approximation of these operators is analysed in $C({\mathbb{R}}).$ Finally, we give some examples of kernels such as B-spline kernels and Blackman-Harris kernel, to which the theory can be applied.

\section{Main Results}

Let $\varphi\in C({\mathbb{R}})$ be fixed. For every $\nu \in {\mathbb{N}}_{0}={\mathbb{N}}\cup\{0\},$ $u\in{\mathbb{R}}$ we define the algebraic moments as
\begin{eqnarray*}
m_{\nu}(\varphi,u):=\sum_{k=-\infty}^{\infty}\varphi(u-k)(k-u)^{\nu}
\end{eqnarray*}
and the absolute moments by
\begin{eqnarray*}
M_{\nu}(\varphi):=\sup_{u\in{\mathbb{R}}}\sum_{k=-\infty}^{\infty}|\varphi(u-k)||(k-u)|^{\nu}.
\end{eqnarray*}

\begin{Remark}
Note that for $\mu, \nu\in\bf{N_0}$ with $\mu<\nu,$ $M_{\nu}(\varphi)<+\infty$ implies $M_{\mu}(\varphi)<+\infty.$ Indeed for $\mu<\nu,$ we have
\begin{eqnarray*}
\sum_{k=-\infty}^{\infty}|\varphi(u-k)||(k-u)|^{\nu}&=&\sum_{|u-k|<1}|\varphi(u-k)||(k-u)|^{\nu}
+\sum_{|u-k|\geq1}|\varphi(u-k)||(k-u)|^{\nu}\\
&\leq&2\|\varphi\|_{\infty}+ \sum_{|u-k|\geq1}|\varphi(u-k)|\frac{|(k-u)|^{\nu}}{|(k-u)|^{{\nu}-\mu}}\\
&\leq&2\|\varphi\|_{\infty}+ M_{\nu}(\varphi).
\end{eqnarray*}
When $\varphi$ is compactly supported, we immediately have that
$M_{\nu}(\varphi)<+\infty$ for every $\nu\in{\mathbb{N}}_0.$
\end{Remark}

We suppose that the following assumptions hold:
\begin{itemize}
\item[(i)] for every $u\in{\mathbb{R}},$ we have
\begin{eqnarray*}
\sum_{k=-\infty}^{\infty}\varphi(u-k)=1,
\end{eqnarray*}
\item[(ii)] $M_2(\varphi)<+\infty$ and there holds
\begin{eqnarray*}
\lim_{r\rightarrow\infty}\sum_{|u-k|>r}|\varphi(u-k)|(k-u)^2=0
\end{eqnarray*}
uniformly with respect to  $u\in{\mathbb{R}},$

\item[(iii)] for every  $u\in{\mathbb{R}},$ we have
\begin{eqnarray*}
m_{1}(\varphi,u)=m_1(\varphi)=\sum_{k=-\infty}^{\infty}\varphi(u-k)(k-u)=0,
\end{eqnarray*}
\item[(iv)] $\displaystyle{\sup_{k}\{|a_{k}|,|b_{k}|\}=M^{*}<\infty}.$
\end{itemize}

\begin{Theorem}\label{t1}
Let $f\in C({\mathbb{R}})$ and $\{a_k\}$ and $\{b_k\}$ be two bounded sequences of real numbers such that $a_k+b_k=\alpha$ and $b_{k}-a_{k}> \Delta^{*}>0.$ If $f'(x)$ exits at a point $x\in{\mathbb{R}}$ then,
\begin{eqnarray*}
\lim_{w\rightarrow\infty}w[(K_w^{\varphi}f)(x)-f(x)]=\frac{\alpha}{2}f'(x).
\end{eqnarray*}

\noindent\bf{Proof.}\rm \
Let $M^{*}=\sup_{k}\{|a_{k}|,|b_{k}|\}.$
From the Taylor's theorem, we have
\begin{eqnarray*}
f(u)=f(x)+f'(x)(u-x)+h(u-x)(u-x),
\end{eqnarray*}
for some bounded function $h$ such that $h(t)\rightarrow 0$ as $t\rightarrow 0.$\\

Thus, we have
\begin{eqnarray*}
(K_w^{\varphi}f)(x)-f(x)&=&f'(x)\sum_{k=-\infty}^{\infty}\frac{w}{\Delta_{k}}\varphi(wx-k)
\int_{\frac{k+a_k}{w}}^{\frac{k+b_k}{w}}(u-x)du\\&&\nonumber+
\sum_{k=-\infty}^{\infty}\frac{w}{\Delta_k}\varphi(wx-k)\int_{\frac{k+a_k}{w}}^{\frac{k+b_k}{w}}h(u-x)(u-x)du\\
&=&I_1+I_2, (\mbox{say}).
\end{eqnarray*}
First, we obtain $I_1.$
\begin{eqnarray*}
I_1&=&f'(x)\sum_{k=-\infty}^{\infty}\frac{w}{b_k-a_k}\varphi(wx-k)\int_{\frac{k+a_k}{w}}^{\frac{k+b_k}{w}}(u-x)du\\
&=&\frac{f'(x)}{2}\sum_{k=-\infty}^{\infty}\frac{w}{b_k-a_k}\varphi(wx-k)
\bigg[\bigg(\frac{k+b_k}{w}-x\bigg)^2-\bigg(\frac{k+a_k}{w}-x\bigg)^2\bigg]\\
&=&\frac{f'(x)}{2w}\sum_{k=-\infty}^{\infty}\varphi(wx-k)[(b_k+a_k)+2(k-wx)]\\
&=&\frac{\alpha f'(x)}{2w}\sum_{k=-\infty}^{\infty}\varphi(wx-k)+\frac{f'(x)}{w}\sum_{k=-\infty}^{\infty}\varphi(wx-k)(k-wx)\\
&=&\frac{\alpha f'(x)}{2w}.
\end{eqnarray*}
Next, we obtain $I_2.$
\begin{eqnarray*}
I_2=\sum_{k=-\infty}^{\infty}\frac{w}{b_k-a_k}\varphi(wx-k)\int_{\frac{k+a_k}{w}}^{\frac{k+b_k}{w}}h(u-x)(u-x)du.
\end{eqnarray*}
In order to  obtain an estimate of $I_2,$ let $\epsilon>0$ be fixed. Then, there exists $\delta>0$ such that $|h(t)|\leq \epsilon$ for  $|t|\leq\delta. $  Then, we have
\begin{eqnarray*}
|I_2|&\leq&\sum_{|k-wx|<\delta w/2}\frac{w}{b_k-a_k}|\varphi(wx-k)|\int_{\frac{k+a_k}{w}}^{\frac{k+b_k}{w}}|h(u-x)||(u-x)|du\\&&+\sum_{|k-wx|\geq\delta w/2}\frac{w}{b_k-a_k}|\varphi(wx-k)|\int_{\frac{k+a_k}{w}}^{\frac{k+b_k}{w}}|h(u-x)||(u-x)|du=J_1+J_2, \mbox{(say)}.
\end{eqnarray*}
First, we estimate $J_1.$

We have
\begin{eqnarray*}
J_1&\leq&{\epsilon}\sum_{|k-wx|<\delta w/2}\frac{w}{b_k-a_k}|\varphi(wx-k)|\int_{\frac{k+a_k}{w}}^{\frac{k+b_k}{w}}|(u-x)|du\\
&\leq&\frac{\epsilon}{2}\sum_{k=-\infty}^{\infty}\frac{w}{b_k-a_k}|\varphi(wx-k)|
\bigg[\bigg(\frac{k+b_k}{w}-x\bigg)^2+\bigg(\frac{k+a_k}{w}-x\bigg)^2\bigg]\\
&\leq&\frac{\epsilon}{w\Delta^{*}}\bigg((M^{*})^{2}M_0(\varphi)+2M^{*} M_1(\varphi)+M_2(\varphi)\bigg).
\end{eqnarray*}
Next, we obtain $J_2.$ Let $R>0$ be such that
\begin{eqnarray*}
\sum_{|u-k|>R}|\varphi(u-k)|(u-k)^2<\epsilon
\end{eqnarray*}
uniformly with respect to $u\in{\mathbb{R}}.$ Also, let ${w}$ be such that $\delta w/2>R.$ Then
\begin{eqnarray*}
\sum_{|k-wx|>\delta w/2}|\varphi(wx-k)|(wx-k)^2<\epsilon
\end{eqnarray*}
for every $x\in {\mathbb{R}}.$ The same inequality holds also for the series
\begin{eqnarray*}
\sum_{|k-wx|>\delta w/2}|\varphi(wx-k)||wx-k|^j<\epsilon
\end{eqnarray*}
for $j=0,1.$ Hence, we get
\begin{eqnarray*}
J_2&\leq&{\|h\|_{\infty}}\sum_{|k-wx|\geq\delta w/2}\frac{w}{b_k-a_k}|\varphi(wx-k)|\int_{\frac{k+a_k}{w}}^{\frac{k+b_k}{w}}|(u-x)|du\\
&\leq&\frac{\|h\|_{\infty}}{2}\sum_{|k-wx|\geq\delta w/2}\frac{w}{b_k-a_k}|\varphi(wx-k)|
\bigg[\bigg(\frac{k+b_k}{w}-x\bigg)^2+\bigg(\frac{k+a_k}{w}-x\bigg)^2\bigg]\\
&\leq&\frac{\epsilon\|h\|_{\infty}}{w\Delta^{*}}(1+M^{*})^2.
\end{eqnarray*}
Hence, the proof is completed.
\end{Theorem}

\begin{Remark}
The boundedness assumption on $f$ can be relaxed by assuming that there are two positive constant $a,b$ such that $|f(x)|\leq a+b|x|,$ for every $x\in {\mathbb{R}}.$\\

We have
\begin{eqnarray*}
|K_w^{\varphi}f)(x)|&\leq&\sum_{k=-\infty}^{\infty}\frac{w}{b_k-a_k}|\varphi(wx-k)|\int_{\frac{k+a_k}{w}}^{\frac{k+b_k}{w}}|f(u)|du\\
&\leq&\sum_{k=-\infty}^{\infty}\frac{w}{b_k-a_k}|\varphi(wx-k)|\int_{\frac{k+a_k}{w}}^{\frac{k+b_k}{w}}(a+b|u|)du\\
&\leq&M_0(\varphi)(a+b|x|) + \frac{b}{\Delta^{*}w}\bigg( M_2(\varphi)+ 2(M^{*}) M_1(\varphi) + (M^{*})^2 M_0(\varphi)\bigg)
\end{eqnarray*}
and hence the series $K_w^{\varphi}f$ is absolutely convergent for every $x\in{\mathbb{R}}.$ Moreover, for a fixed $x_0\in{\mathbb{R}},$
\begin{eqnarray*}
P_1(x)=f(x_0)+f'(x_0)(x-x_0),
\end{eqnarray*}
the Taylor's polynomial of first order centered at the point $x_0,$ by the Taylor's formula we can write
\begin{eqnarray*}
\frac{f(x)-P_1(x)}{(x-x_0)}=h(x-x_0),
\end{eqnarray*}
where $h$ is a function such that $\displaystyle\lim_{t\rightarrow 0}h(t)=0.$ Then $h$ is bounded on $[x_0-\delta, x_0+\delta],$ for some $\delta>0.$ For $|x-x_0|>\delta,$ we have
\begin{eqnarray*}
|h(x-x_0)|\leq \frac{a+b|x|}{|x-x_0|}+\frac{|P_1(x)|}{|x-x_0|}\le \frac{a+b|x|}{|x-x_0|}+\frac{|f(x_{0})|}{|x-x_{0}|}+|f'(x_{0})|,
\end{eqnarray*}
and the  terms on the right-hand side of the above inequality are all bounded for $|x-x_0|>\delta.$ Hence, $h(.-x_0)$ is bounded on ${\mathbb{R}}.$ Along the lines of the proof of Theorem \ref{t1}, the same Voronovskaya formula can be obtained.
\end{Remark}

\subsection{Quantitative Estimate}
Let $C^{m}$ denote the set of all $f\in C({\mathbb{R}})$ such that $f$ is $m$ times continuously differentiable and $\|f^{(m)}\|_{\infty}<\infty.$

Let $\delta>0.$ For $f\in C({\mathbb{R}}),$ the Peetre's $K$-functional is defined as
\begin{eqnarray*}
K(\delta,f,C,C^{1}):=\inf\{\|f-g\|_{\infty}+\delta\|g'\|_{\infty}: g\in C^{1}\}.
\end{eqnarray*}
For a given $\delta>0,$ the usual modulus of continuity of a given uniformly continuous function $f:{\mathbb{R}}\rightarrow{\mathbb{R}}$ is defined as
\begin{eqnarray*}
\omega(f,\delta):=\sup_{|x-y|\leq \delta}|f(x)-f(y)|.
\end{eqnarray*}
It is well known that, for any positive constant $\lambda>0$, the modulus of continuity satisfies the following property
\begin{eqnarray}\label{e2}
\omega(f,\lambda\delta)\leq (\lambda+1)\omega(f,\delta).
\end{eqnarray}
For a function $f\in C^{m},$ $x_0,x\in{\mathbb{R}}$ and $m\geq 1$, the Taylor's formula is given by
\begin{eqnarray*}
f(x)=\sum_{k=0}^{m}\frac{f^{(k)}(x_0)}{k!}(x-x_0)^k+R_m(f;x_0,x)
\end{eqnarray*}
and the remainder term $R_m(f;x_0,x)$ is estimated by
\begin{eqnarray*}
|R_m(f;x_0,x)|\leq\frac{|x-x_0|^m}{m!}\omega(f^{(m)};|x-x_0|).
\end{eqnarray*}
For every $f\in C({\mathbb{R}})$ there holds
\begin{eqnarray}\label{e3}
K(\delta/2,f,C,C^{1})=\frac{1}{2}\overline{\omega}(f,\delta),
\end{eqnarray}
where $\overline{\omega}(f,.)$ denotes the least concave majorant of ${\omega}(f,.)$ (see e.g. \cite{ana}).

The following estimate for the remainder $R_m(f;x_0,x)$ in terms of $\overline{\omega}$ was proved in \cite{gonska}.

\begin{Lemma}\label{l1}
Let $f\in C^{m},$ $m\in\mathbb{N}^{0}$ and $x_0,x\in{\mathbb{R}}.$ Then, we have
\begin{eqnarray*}
|R_m(f;x_0,x)|\leq\frac{|x-x_0|^m}{m!}\overline{\omega}\bigg(f^{(m)};\frac{|x-x_0|}{m+1}\bigg).
\end{eqnarray*}
\end{Lemma}
We have the following quantitative version of Theorem \ref{t1} in terms of the modulus $\overline{\omega},$ in case of $m=1.$

\begin{Theorem}\label{t2}
Let $f\in C^{1}$ and $\{a_k\}$ and $\{b_k\}$ be two sequences of real numbers such that $a_{k}+b_{k}=\alpha,$ $b_{k}-a_{k}\ge \Delta^{*}>0$  and $\sup_{k}\left\{|a_{k}|,|b_{k}|\right\}\le M^{*}$. Then, for very $x\in{\mathbb{R}},$ the following hold: \begin{eqnarray*}
\left|w[(K_w^{\varphi}f)(x)-f(x)]-\frac{\alpha f'(x)}{2}\right|\leq \frac{A}{\Delta^{*}}\overline{\omega}\bigg(f',\frac{\Delta^{*}}{2w}\bigg)
\end{eqnarray*}
where $A=(M^{*})^2 M_0(\varphi)+2 (M^{*})M_1(\varphi)+ M_2(\varphi).$

\noindent\bf{Proof.}\rm \
Let $f\in C^{1}$ be fixed. Then, we can write\\

\noindent
$\left|w[(K_w^{\varphi}f)(x)-f(x)]-\frac{\alpha f'(x)}{2}\right|$
\begin{eqnarray*}
&=&\bigg|{f'(x)}\sum_{k=-\infty}^{\infty}\frac{w^2}{b_k-a_k}\varphi(wx-k)\int_{\frac{k+a_k}{w}}^{\frac{k+b_k}{w}}(u-x)du\\&&+
\sum_{k=-\infty}^{\infty}\frac{w^2}{b_k-a_k}\varphi(wx-k)\int_{\frac{k+a_k}{w}}^{\frac{k+b_k}{w}}h(u-x)(u-x)du-\frac{\alpha f'(x)}{2}\bigg|\\
&\leq&\sum_{k=-\infty}^{\infty}\frac{w^2}{b_k-a_k}|\varphi(wx-k)|\int_{\frac{k+a_k}{w}}^{\frac{k+b_k}{w}}|h(u-x)||(u-x)|du.
\end{eqnarray*}
Using the relation (\ref{e3}) and Lemma \ref{l1}, we obtain\\

\noindent
$\left|w[(K_w^{\varphi}f)(x)-f(x)]-\frac{\alpha f'(x)}{2}\right|$
\begin{eqnarray*}
&\leq&\sum_{k=-\infty}^{\infty}\frac{w^2}{b_k-a_k}|\varphi(wx-k)|\int_{\frac{k+a_k}{w}}^{\frac{k+b_k}{w}}|(u-x)|\overline{\omega}\bigg(f',\frac{|x-u|}{2}\bigg)du\\
&=&2\sum_{k=-\infty}^{\infty}\frac{w^2}{b_k-a_k}|\varphi(wx-k)|\int_{\frac{k+a_k}{w}}^{\frac{k+b_k}{w}}|(u-x)|K\bigg(\frac{|u-x|}{4},f',C,C^{1}\bigg)du:=I_1.
\end{eqnarray*}
For $g\in C^2,$ we have
\begin{eqnarray*}
I_1&\leq& \sum_{k=-\infty}^{\infty}\frac{2w^2}{b_k-a_k}|\varphi(wx-k)|\int_{\frac{k+a_k}{w}}^{\frac{k+b_k}{w}}|(u-x)|
\bigg(\|(f-g)'\|_{\infty}+\frac{|u-x|}{4}\|g''\|_{\infty}\bigg)du\\
&\leq& \|(f-g)\|_{\infty}\sum_{k=-\infty}^{\infty}\frac{2w^2}{b_k-a_k}|\varphi(wx-k)|
\int_{\frac{k+a_k}{w}}^{\frac{k+b_k}{w}}|(u-x)|du\\
&&+
\|g''\|_{\infty}\sum_{k=-\infty}^{\infty}\frac{2w^2}{b_k-a_k}|\varphi(wx-k)|\int_{\frac{k+a_k}{w}}^{\frac{k+b_k}{w}}(u-x)^2du\\
&\leq&\|(f-g)'\|_{\infty}\sum_{k=-\infty}^{\infty}|\varphi(wx-k)|\frac{w^2}{b_k-a_k}
\bigg[\bigg(\frac{k+b_k}{w}-x\bigg)^2+\bigg(\frac{k+a_k}{w}-x\bigg)^2\bigg]\\&&+
\|g''\|_{\infty}\sum_{k=-\infty}^{\infty}\frac{w^2}{6(b_k-a_k)}|\varphi(wx-k)|\bigg[\bigg(\frac{k+b_k}{w}-x\bigg)^3-\bigg(\frac{k+a_k}{w}-x\bigg)^3\bigg]\\
&\leq&\frac{\|(f-g)'\|_{\infty}}{\Delta^{*}}\bigg(2 (M^{*})^2 M_0(\varphi)+4M_1(\varphi)(M^{*})+2M_2(\varphi)\bigg)\\&&+
\|g''\|_{\infty}\frac{1}{6w}\bigg(3 (M^{*})^2 M_0(\varphi)+6 (M^{*})M_1(\varphi)+3M_2(\varphi)\bigg)\\
&\leq&\|(f-g)'\|_{\infty} \frac{2}{\Delta^{*}}\bigg( (M^{*})^2 M_0(\varphi)+2M_1(\varphi)(M^{*})+M_2(\varphi)\bigg) \\&&+
\|g''\|_{\infty}\frac{1}{2 w}\bigg( (M^{*})^2 M_0(\varphi)+2M_1(\varphi)(M^{*})+M_2(\varphi)\bigg)\\
&\leq&\frac{2A}{\Delta^{*}}\bigg(\|(f-g)'\|_{\infty}+\|g''\|_{\infty}\frac{\Delta^{*}}{4w}\bigg).
\end{eqnarray*}
Taking the infimum over all $g\in C^2,$ we get
\begin{eqnarray}
I_1&\leq&\frac{A}{\Delta^{*}}\overline{\omega}\bigg(f',\frac{\Delta^{*}}{2w}\bigg).
\end{eqnarray}
Hence, the proof is completed.
\end{Theorem}

\begin{Remark}
As a consequence of Theorem \ref{t2}, under the above assumptions we get the uniform convergence for
$w[(K_w^{\varphi}f)(x)-f(x)]$ to $\frac{\alpha}{2} f'(x).$
\end{Remark}

\begin{Remark}
Note that when $\varphi$ is supported in $I=[-R,R],$ $R>0$ we can obtain a different estimate for $I_1.$
\begin{eqnarray*}
\left|w[(K_w^{\varphi}f)(x)-f(x)]-\frac{\alpha} {2}f'(x)\right|
&\leq& \frac{M_0(\varphi)(R^2+ 2 R M^{*}+ (M^{*})^2)}{\Delta^{*}}\overline{\omega}\bigg(f',\frac{\Delta^{*}}{2w}\bigg).
\end{eqnarray*}
Also, we obtain
\begin{eqnarray*}
I_1&\leq& \frac{2 M_0(\varphi)(R^2+ 2 R M^{*}+ (M^{*})^2)}{\Delta^{*}}\left[ \|(f-g)'\|_{\infty}+ \|g''\|_{\infty}\frac{\Delta^{*}}{4 w}\right].
\end{eqnarray*}
\end{Remark}

\subsection{Order of approximation}
The order of approximation for the generalized sampling Kantorovich series has been extensively studied by many authors (see \cite{mant}, \cite{Dani}, \cite{vinti1}, \cite{vinti}).

\begin{Theorem}\label{t3}
Let $\varphi$ be a kernel satisfying an additional  condition that $M_{\beta}(\varphi)=\sup_{u\in{\mathbb{R}}}\sum_{k\in{\mathbb{Z}}}|\varphi(u-k)||(k-u)|^{\beta}<+\infty $ for some  $0<\beta<1$ and $\{a_k\}$ and $\{b_k\}$ be two bounded sequences of real numbers. Then, for any $f\in C({\mathbb{R}}),$ we have
\begin{eqnarray*}
|(K_w^{\varphi}f)(x)-f(x)|\leq\omega(f,w^{-\beta})\bigg(M_\beta(\varphi)+(M^{*})^{\beta}M_0(\varphi)+
M_0(\varphi)\bigg)+2^{\beta+1}\|f\|_{\infty}w^{-\beta}M_\beta(\varphi),
\end{eqnarray*}
for every $x\in {\mathbb{R}}$ and $w>2M^{*}.$

\noindent\bf{Proof.}\rm \
Let $x\in{\mathbb{R}}$ be fixed. Then, for $w>0,$ we can write
\begin{eqnarray*}
|(K_w^{\varphi}f)(x)-f(x)|&=&\bigg|(K_w^{\varphi}f)(x)-f(x)\sum_{k\in\bf{Z}}\varphi(wx-k)\bigg|\\
&\leq& \sum_{k\in\bf{Z}}\bigg(\frac{w}{b_k-a_k}\int_{\frac{k+a_k}{w}}^{\frac{k+b_k}{w}}|f(u)-f(x)|du\bigg)
|\varphi(wx-k)|=:J.
\end{eqnarray*}
Now, we estimate $J.$
\begin{eqnarray*}
J&\leq& \sum_{|wx-k|\leq w/2}\bigg(\frac{w}{b_k-a_k}\int_{\frac{k+a_k}{w}}^{\frac{k+b_k}{w}}|f(u)-f(x)|du\bigg)
|\varphi(wx-k)|\\&&+
\sum_{|wx-k|>w/2}\bigg(\frac{w}{b_k-a_k}\int_{\frac{k+a_k}{w}}^{\frac{k+b_k}{w}}|f(u)-f(x)|du\bigg)
|\varphi(wx-k)|=:I_1+I_2.
\end{eqnarray*}
We observe that, for every $u\in\bigg[\frac{k+a_k}{w}, \frac{k+b_k}{w}\bigg]$ and  suitable large $w$ with $|wx-k|\leq w/2,$
we get
\begin{eqnarray*}
|u-x|\leq\left|u-\bigg(\frac{k+a_k}{w}\bigg)+\bigg(\frac{k+a_k}{w}\bigg)-x\right|
\leq\bigg|u-\frac{k}{w}\bigg|+\bigg|\frac{k}{w}-x\bigg|\leq\frac{M^{*}}{w}+\frac{1}{2}\leq1,
\end{eqnarray*}
for every $w\geq 2M^{*}.$ Since $0<\beta<1,$ we have
\begin{eqnarray*}
\omega(f,|u-x|)\leq \omega(f,|u-x|^{\beta}).
\end{eqnarray*}
Using the property of modulus of continuity \ref{e2}, we obtain
\begin{eqnarray*}
I_1&\leq&\sum_{|wx-k|\leq w/2}\bigg(\frac{w}{b_k-a_k}\int_{\frac{k+a_k}{w}}^{\frac{k+b_k}{w}}\omega(f,|u-x|^{\beta})du\bigg)|\varphi(wx-k)|\\
&\leq&\sum_{|wx-k|\leq w/2}\bigg(\frac{w}{b_k-a_k}\int_{\frac{k+a_k}{w}}^{\frac{k+b_k}{w}}
(1+w^{\beta}|u-x|^{\beta})\omega(f,w^{-\beta})du\bigg)|\varphi(wx-k)|\\
&\leq&\omega(f,w^{-\beta})\bigg[\sum_{|wx-k|\leq w/2}\bigg(\frac{w}{b_k-a_k}\int_{\frac{k+a_k}{w}}^{\frac{k+b_k}{w}}w^{\beta}|u-x|^{\beta}du\bigg)|\varphi(wx-k)|+\sum_{|wx-k|\leq/2}|\varphi(wx-k)|\bigg]\\
&=:&\omega(f,w^{-\beta})(J_1+J_2).
\end{eqnarray*}
First, we obtain $J_1.$ Using the property of sub-addivity of $|.|^{\beta}$ with
$0<\beta<1,$ we have
\begin{eqnarray*}
J_1&\leq&\sum_{|wx-k|\leq w/2}\bigg(w^{\beta}\max_{u\in\left[\frac{k+a_k}{w},\frac{k+b_k}{w}\right]}|u-x|^{\beta}\bigg)|\varphi(wx-k)|\\
&\leq&\sum_{|wx-k|\leq w/2}\bigg(w^{\beta}\max\bigg(\left|\frac{k+a_k}{w}-x\right|^{\beta}, \left|\frac{k+b_k}{w}-x\right|^{\beta}\bigg)\bigg)|\varphi(wx-k)|\\
&\leq&\sum_{|wx-k|\leq w/2}\bigg(w^{\beta}\max\bigg(\left|\frac{k}{w}-x\right|^{\beta}+\left|\frac{a_k}{w}\right|^{\beta}, \left|\frac{k}{w}-x\right|^{\beta}+\left|\frac{b_k}{w}\right|^{\beta}\bigg)\bigg)|\varphi(wx-k)|\\
&\leq&\sum_{|wx-k|\leq w/2}w^{\beta}\bigg(\left|\frac{k}{w}-x\right|^{\beta}+\left(\sup_{k}\{|a_k|,|b_{k}|\}\right)^{\beta}w^{-\beta}\bigg)|\varphi(wx-k)|\\
&\leq&\sum_{|wx-k|\leq w/2}|k-wx|^{\beta}|\varphi(wx-k)|+\sum_{|wx-k|\leq w/2}(M^{*})^{\beta}
|\varphi(wx-k)|\\
&\leq&M_{\beta}(\varphi)+(M^{*})^{\beta}M_{0}(\varphi)<\infty.
\end{eqnarray*}
It is easy to see that
\begin{eqnarray*}
J_2&\leq&\sum_{|wx-k|\leq w/2}|\varphi(wx-k)|=M_{0}(\varphi).
\end{eqnarray*}
Next, we estimate $I_2.$
\begin{eqnarray*}
I_2&\leq&2\|f\|_{\infty}\sum_{|wx-k|>w/2}\bigg(\frac{w}{b_k-a_k}\int_{\frac{k+a_k}{w}}^{\frac{k+b_k}{w}}du\bigg)|\varphi(wx-k)|\\
&\leq&2\|f\|_{\infty}\sum_{|wx-k|>w/2}|\varphi(wx-k)|\\
&\leq&2\|f\|_{\infty}\sum_{|wx-k|>w/2}\frac{|wx-k|^{\beta}}{|wx-k|^{\beta}}|\varphi(wx-k)|\\
&\leq&\frac{2\|f\|_{\infty}}{w^{\beta}}\sum_{|wx-k|>w/2}{|wx-k|^{\beta}}|\varphi(wx-k)|\\
&\leq&2^{\beta+1}\|f\|_{\infty}w^{-\beta}M_{\beta}(\varphi)<+\infty,
\end{eqnarray*}
which completes the proof.
\end{Theorem}

\subsection{Applications to special kernels.}
In this section, we describe some particular examples of kernels $\varphi$ which illustrates the previous theory. In particular, we will examine the $B$-splines kernel and Blackman-Harries kernel.
\subsection{Combinations of B-Spline Functions }
First, we consider the sampling Kantorovich operators based upon the combinations of spline functions. For $h\in\mathbb{N},$ the $B-$spline of order $h$ is defined as
\begin{eqnarray*}
B_h(x):=\chi_{[-\frac{1}{2}, \frac{1}{2}]}\star \chi_{[-\frac{1}{2}, \frac{1}{2}]}\star \chi_{[-\frac{1}{2}, \frac{1}{2}]}\star...\star\chi_{[-\frac{1}{2}, \frac{1}{2}]},\mbox{( h \,\ \mbox{times})}
\end{eqnarray*}
where
\begin{eqnarray*}\chi_{[-\frac{1}{2}, \frac{1}{2}]}=
\left\{
\begin{array}{ll}
1 , &\mbox{if} -\frac{1}{2}\leq x\leq\frac{1}{2}\\
0, & \mbox{otherwise}
\end{array}
\right.
\end{eqnarray*}
and $*$ denotes the convolution.\\
The Fourier transform of the functions $B_h(x)$ is given by
\begin{eqnarray*}
\widehat{B}_h(w)=\bigg(\widehat{\chi_{[\frac{-1}{2}, \frac{1}{2}]}}(w)\bigg)^{h}
=\bigg(\frac{\sin w/2}{w/2}\bigg)^{h}, \,\ w\in{\mathbb{R}}, \,\ h\in\mathbb{N}
\end{eqnarray*}
(see \cite{Stens2} and \cite{scho}). Given real numbers $\epsilon_0, \epsilon_1$ with
$\epsilon_0<\epsilon_1$ we will construct the linear combination of translates of $B_h$, with $h\geq 2$ of type
\begin{eqnarray*}
\varphi(x)=a_0B_h(x-\epsilon_0)+a_1B_h(x-\epsilon_1).
\end{eqnarray*}
The Fourier transform of $\varphi$ is given by
\begin{eqnarray*}
\widehat{\varphi}(w)=\bigg(a_0e^{-i\epsilon_0 w}+a_1e^{-i\epsilon_1 w}\bigg)\widehat{B}_{h}(w).
\end{eqnarray*}
Using the Poisson summation formula
\begin{eqnarray*}
(-i)^{j}\sum_{k=-\infty}^{\infty}\varphi(u-k)(u-k)^{j}\sim \sum_{k=-\infty}^{\infty}\widehat{\varphi}^{(j)}(2\pi k)e^{i2\pi ku},
\end{eqnarray*}
we obtain
\begin{eqnarray*}
\sum_{k=-\infty}^{\infty}\varphi(u-k)=\sum_{k=-\infty}^{\infty}\widehat{\varphi}(2\pi k)e^{i2\pi ku}.
\end{eqnarray*}
We have
\begin{eqnarray*}
\widehat{B_{h}}(2\pi k)=\bigg(\frac{\sin(\pi k)}{\pi k}\bigg)^{h}
=
\left\{
\begin{array}{ll}
1 , &\mbox{if}\,\ k=0\\
0, & \mbox{if}\,\ k\neq 0
\end{array}
\right.
\end{eqnarray*}
and hence
\begin{eqnarray*}
\widehat{\varphi}(2\pi k)
=
\left\{
\begin{array}{ll}
a_0+a_1 , &\mbox{if}\,\ k=0\\
0, & \mbox{if}\,\ k\neq 0.
\end{array}
\right.
\end{eqnarray*}
Thus
\begin{eqnarray*}
\sum_{k=-\infty}^{\infty}\varphi(u-k)=a_0+a_1.
\end{eqnarray*}
Therefore, condition (i) is satisfied if $a_0+a_1=1.$ Now, we show that condition (iii) is also satisfied.

Again from the Poisson summation formula, we obtain
\begin{eqnarray*}
(-i)\sum_{k=-\infty}^{\infty}\varphi(u-k)(u-k)=
\sum_{k=-\infty}^{\infty}\widehat{\varphi}^{\prime}(2\pi k)e^{i2\pi ku}.
\end{eqnarray*}
Also, we have
\begin{eqnarray*}
\widehat{\varphi}^{\prime}(w)=(-i\epsilon_0a_0e^{-i\epsilon_0 w}-i\epsilon_1a_1e^{-i\epsilon_1 w})\widehat{B}_{h}(w)+(a_0e^{-i\epsilon_0 w}+a_1e^{-i\epsilon_1 w})\widehat{B}_{h}^{\prime}(w)
\end{eqnarray*}
Since $\widehat{B}_{h}^{\prime}(2\pi k)=0,$ $\forall k$ which implies that $\widehat{\varphi}^{\prime}(2\pi k)=0.$ Thus, we have
\begin{eqnarray*}
\widehat{\varphi}(0)=a_0+a_1=1,\,\,\ \widehat{\varphi}^{\prime}(0)=\epsilon_0a_0+\epsilon_1a_1=0.
\end{eqnarray*}
Solving the above linear system we get the unique solution
\begin{eqnarray*}
a_0=\frac{\epsilon_1}{\epsilon_1-\epsilon_0}, a_1=-\frac{\epsilon_0}{\epsilon_0-\epsilon_1}.
\end{eqnarray*}
Moreover it is easy to see that the support of the function $\varphi$ is contained in the interval $[\epsilon_0-\frac{h}{2}, \epsilon_1-\frac{h}{2},].$
Since $\varphi(u-k)=0$ if $|u-k|>r$ for $r$ sufficiently large, we have
\begin{eqnarray*}
\lim_{r\rightarrow\infty}\sum_{|k-u|>r}\varphi(u-k)(k-u)^2=0.
\end{eqnarray*}
Condition (ii) is satisfied. Finally, we verify the condition that $M^{\beta}(\varphi)<\infty.$
\begin{eqnarray*}
\sum_{k\in \mathbb{Z}}|\varphi(u-k)||(k-u)|^{\beta}
&=& \sum_{|k-u|<R}|\varphi(u-k)||(k-u)|^{\beta}\\&&+
\sum_{|k-u|\geq R}|\varphi(u-k)||(k-u)|^{\beta}.
\end{eqnarray*}
We can see that $\displaystyle\sup_{u}|{\{k: |u-k|<R\}}|\leq N_0.$ Thus, we get
\begin{eqnarray*}
M_{\beta}(\varphi)=\sum_{k\in \mathbb{Z}}|\varphi(u-k)||(k-u)|^{\beta}<\infty.
\end{eqnarray*}

\subsection{A particular Blackman-Harris kernel}
Next, we consider the Blackman-Harris kernel.  For every $x\in{\mathbb{R}},$ we define the kernel
(see \cite{man})
\begin{eqnarray*}
\varphi(x)\equiv H(x)=\frac{1}{2}sinc(x)+\frac{9}{32}(sinc(x+1)+sinc(x-1))-
\frac{1}{32}(sinc(x+3)+sinc(x-3)),
\end{eqnarray*}
where $sinc(x)=\frac{\sin \pi x}{\pi x}.$ From \cite{a5}, there holds that $H(x)=O(|x|^{-5})$ as $|x|\rightarrow\infty.$ In view of \cite{BVBS}, it follows that $M_2(H)$ is finite and
\begin{eqnarray*}
\lim_{r\rightarrow\infty}\sum_{|k-u|>r}|H(u-k)|(u-k)^2=0.
\end{eqnarray*}
Indeed, there exists $N_0>0$ such that $|H(x)|\leq{M}/|x|^5$ for $|x|\geq N_0.$ Thus, we have for $r>N_0$
\begin{eqnarray*}
\sum_{|k-u|>r}|H(u-k)|(u-k)^2\leq M\sum_{|k-u|>r}\frac{1}{|u-k|^3}\leq \frac{M}{r}
\sum_{|k-u|>r}\frac{1}{|u-k|^2}\leq\frac{2M}{r} \sum_{k=1}^{\infty}\frac{1}{k^2}.
\end{eqnarray*}
The Fourier transform of the function $H(x)$ is given by
\begin{eqnarray*}
\widehat{H}(w)=\frac{1}{\sqrt{2\pi}}\lambda\bigg(\frac{w}{\pi}\bigg),
\end{eqnarray*}
where $\lambda(w)=\bigg(\frac{1}{2}+\frac{9}{16}\cos(\pi w)-\frac{1}{16}\cos(3\pi w)\bigg)\chi_{[-1,1]}(w),$ $\chi_{I}$ is the characteristic function of the set $I.$ From Lemma 3 in \cite{Stens2}, we obtain
\begin{eqnarray*}
m_1{(H)}=\sum_{k=-\infty}^{\infty}H(u-k)(u-k)=0.
\end{eqnarray*}
Hence the condition (i)-(iii) are satisfied. Finally, we verify that $M^{\beta}(\varphi)<\infty.$
\begin{eqnarray*}
\sum_{k\in \mathbb{Z}}|H(u-k)||(k-u)|^{\beta}
&=& \sum_{|k-u|<R}|H(u-k)||(k-u)|^{\beta}\\&&+
\sum_{|k-u|\geq R}|H(u-k)||(k-u)|^{\beta}=S_1+S_2, (\mbox{say}).
\end{eqnarray*}
First, we consider $S_2.$ There exists $N>0$ such that $|H(x)|\leq{M}/|x|^5$ for $|x|\geq N.$ Thus, we have for $R>N,$
\begin{eqnarray*}
S_2\leq M\sum_{|k-u|\geq R}\frac{|u-k|^{\beta}}{|u-k|^5}\leq{2M}\sum_{k=1}^{\infty}\frac{1}{k^{5-\beta}}.
\end{eqnarray*}
Next, we estimate $S_1.$ We have $\displaystyle\sup_{u}|{\{k: |u-k|<R\}}|<\infty.$ Thus, we obtain
\begin{eqnarray*}
S_1\leq \sum_{|k-u|<R}|H(u-k)||k-u|^{\beta}\leq\sum_{|k-u|<R}M_0 R^{\beta}\leq
\lceil 2R\rceil M_0R^{\beta},
\end{eqnarray*}
where $\lceil x\rceil$ denotes the smallest integer greater than or equal to $x.$
Hence, we get
\begin{eqnarray*}
M_{\beta}(H)=\sum_{k\in \mathbb{Z}}|H(u-k)||(k-u)|^{\beta}<\infty.
\end{eqnarray*}
Thus, all the conditions are satisfied for the function $H(x).$

\end{document}